\documentclass{article}%
\usepackage{amsmath}
\usepackage{amsfonts}
\usepackage{amssymb}
\usepackage{graphicx}%
\newtheorem{theorem}{Theorem}

\newtheorem{definition}{Definition}

\newtheorem{lemma}{Lemma}

\newtheorem{proposition}{Proposition}

\begin{document}

\title{On discrimination between two close distribution tails.}
\author{Rodionov~I.\,V.\thanks{Moscow State University, Faculty of Mathematics and Mechanics and Moscow Institute of Physics and Technology, Faculty of Innovations and High Technologies. E-mail: vecsell@gmail.com}}
\maketitle

\section{Introduction. Main result.}

Statistics deals often with discrimination of close distributions based on
censored or truncated data, in particular, for high-risk insurances and
reliability problems. The situation when one observes data exceeding a
pre-determined threshold is well-studied, see \cite{Dufour}, \cite{Guilbaud},
\cite{Chernobai} and references therein. On the other hand statistics of
extremes says that only higher order statistics should be used for
discrimination of close distribution tails, wherein moderate sample values can
be modeled with standard statistical tools. In particular, such approach for
distributions from Gumbel maximum domain of attraction (for the definitions
see \cite{Dehaan}) is considered in \cite{girard}, \cite{Rodionov1},
\cite{Rodionov2}. As well, any estimators of the extreme value indices
$\gamma$ and $\rho$ (see \cite{HaanResnick}) can be used also to discriminate
the distribution tails. Notice that  we do not assume belonging the
corresponding distribution function to a maximum domain of attraction.
\newline

\begin{definition}
\label{def1}The distribution functions $F$ and $G$ are said to be satisfied
the condition $B(F,G)$ if for some $\varepsilon>0$ and $x_{0}$
\begin{equation}
\frac{1-G(x)}{(1-F(x))^{1-\varepsilon}}\text{ is nondecreasing with}\ x>x_{0}.
\label{cond}%
\end{equation}

\end{definition}
Denote by $\Theta(F_{0})$ the class of continuous distribution functions
$F_{1}$ satisfying either $B(F_{1},F_{0})$ or $B(F_{0},F_{1}).$ Consider the
simple hypothesis $H_{0}:F=F_{0}$ and the alternative hypothesis $H_{1}%
:F\in\Theta(F_{0}),$ where $F_{0}$ is continuous. Notice that if distribution functions $F,$ $G$ satisfy either $B(F,G)$ or $B(G,F)$ for some $\varepsilon>0$ then it holds for all $\varepsilon_1,\ 0<\varepsilon_1<\varepsilon.$ So denote
\[\varepsilon(F, G) = \max\{\varepsilon:\ F, G \text{ satisfy either } B(F,G) \text{ or } B(G,F) \text{ for } \varepsilon\}. \]Denote by $\Theta_\varepsilon(F_{0})$ the class of continuous distribution functions
$F_{1}$ satisfying either $B(F_{1},F_{0})$ or $B(F_{0},F_{1})$ with $\varepsilon(F_0, F_1)\geq\varepsilon$ and consider another alternative hypothesis $H_{1, \varepsilon}: F\in\Theta_\varepsilon(F_{0}).$ Let $X_{1},\ldots,X_{n}$ be
i.i.d. random variables with a common distribution function $F$. Denote by
$X_{(1)}\leq\ldots\leq X_{(n)}$ the order statistics for them. Introduce the
Hill-like statistics
\[
R_{k,n}=\ln(1-F_{0}(X_{(n-k)}))-\frac{1}{k}\sum\limits_{i=n-k+1}^{n}%
\ln(1-F_{0}(X_{(i)})).
\]
which we are going to use for the problem of discrimination between the two
introduced above hypotheses when $k$ higher order statistics are known. Remark
that if $F_{0}$ is Pareto distribution function with parameter $\gamma$, then
\[
R_{k,n}\overset{d}{=}\gamma_{H}/\gamma,
\]
where $\gamma_{H}$ is the Hill estimator of $\gamma.$ If furthermore $F_{0}$
belongs to Fr\'{e}chet max-domain of attraction, then $R_{k,n}$ behaves
asymptotically as $\gamma_{H}/\gamma,$ that is, theirs ratio tends to one as
$n\rightarrow\infty.$ We will show that the distributions of $R_{k,n}$ if
either $H_{0}$ or $H_{1}$ fulfilled are different which can give a statistical
for discrimination the hypotheses. The following two results describe the
behavior of $R_{k,n}$ as $k,n\rightarrow\infty$ with $k<n$ provided $H_{0}$ or
$H_{1}$ is fulfilled.

\begin{theorem}
\label{Th1} If $H_{0}$ holds then
\[
\sqrt{k} (R_{k,n} - 1) \overset{d}{\longrightarrow} \xi\text{ as } k,
n\to\infty,
\]
where $\xi$ is standard normal random variable, i.e. $\xi\sim N(0,1).$
\end{theorem}

This theorem gives obvious goodness-of-fit test for the tail of $F.$ Besides,
the following result provides some information about the consistency of this
test. Assume that $H_{0}$ does not hold and $F$ is equal to $F_{1}$ which is
different from $F_{0}.$ Denote $x^{\ast}$, the right endpoint of $F_{1}$, that
is, $x^{\ast}=\sup\{x:F_{1}(x)<1\}.$ Assume that $F_{0}$ and any $F_{1}%
\in\Theta(F_{0})$ have the same right endpoint (how to discriminate
distributions with different endpoints, see \cite{Falk}, \cite{Dehaan}).
Further consider $x^{\ast}=+\infty,$ otherwise change variables $y=1/(x^{\ast
}-x)$ gives the assumption. The following theorem shows consistency of the
proposed test.

\begin{theorem}
\label{Th2}
\begin{enumerate}
\item[\bf (i)] If $H_{1}$ holds then
\[
\sqrt{k_{n}} |R_{k_{n},n} - 1| \overset{d}{\longrightarrow} +\infty
\]
provided $k_{n}\to\infty,$ $k_{n}/n\to0$ as $n\to\infty.$

\item[\bf (ii)] If $H_{1,\varepsilon}$ holds then under the same conditions
\[
\inf\limits_{F_1\in\Theta_{\varepsilon}(F_0)}\sqrt{k_{n}} |R_{k_{n},n} - 1| \overset{d}{\longrightarrow} +\infty.
\]
\end{enumerate}
\end{theorem}

The considered test makes it possible to discriminate, for example, two normal
distributions with different variances, but we should weaken the condition
\eqref{cond} to discriminate two normal distributions with the same variance
and different means. But weakening the condition \eqref{cond} imposes some
conditions on behavior of the sequence $k_{n}.$

\begin{definition}
The distribution functions $F$ and $G$ are said to satisfy the condition $C(F,
G)$ if for some $\varepsilon>0$ and $x_{0}$
\begin{equation}
\frac{1-G(x)}{(1-F(x))(-\ln(1-F(x)))^{\varepsilon}} \text{ is nondecreasing,
}\ \ x>x_{0}. \label{cond1}%
\end{equation}

\end{definition}
Denote by $\Theta^{\prime}(F_{0})$ the class of continuous distribution functions
$F_{1}$ satisfying either $C(F_{1}, F_{0})$ or $C(F_{0}, F_{1})$ and the following
condition: for some $\delta\in(0,1)$
\begin{equation}
1 - F_{1}(x)\leq(1-F_{0}(x))^{\delta},\ \ x>x_{0}. \label{new}%
\end{equation}
See, if distribution functions $F,$ $G$ satisfy either $C(F,G)$ or $C(G,F)$ for some $\varepsilon>0$ then it holds for all $\varepsilon_1,\ 0<\varepsilon_1<\varepsilon.$ Denote
\[\varepsilon^{\prime}(F, G) = \max\{\varepsilon:\ F, G \text{ satisfy either } C(F,G) \text{ or } C(G,F) \text{ for } \varepsilon\}. \]Denote by $\Theta_\varepsilon^{\prime}(F_{0})$ the class of continuous distribution functions
$F_{1}$ satisfying \eqref{new} and either $C(F_{1},F_{0})$ or $C(F_{0},F_{1})$ with $\varepsilon^{\prime}(F_0, F_1)\geq\varepsilon.$
As before, consider the simple hypothesis $H_{0}: F = F_{0}$ and two
alternative hypotheses $H_{1}^{\prime}: F \in\Theta^{\prime}(F_{0}),$ $H_{1,\varepsilon}^{\prime}: F \in\Theta^{\prime}_\varepsilon(F_{0})$ with
continuous $F_{0}.$

\begin{theorem}
\label{Th3}
\begin{enumerate}
\item[\bf (i)] If $H_{1}^{\prime}$ holds then
\[
\sqrt{k_{n}} |R_{k_{n},n} - 1| \overset{d}{\longrightarrow} +\infty
\]
provided $k_{n}/n\to0,$ $k_{n}^{1/2-\alpha}/\ln n \to+\infty,$ for some
$\alpha\in(0, 1/2),$ as $n\to\infty.$
\item[\bf (ii)] If $H_{1,\varepsilon}^{\prime}$ holds then under the same conditions
\[
\inf_{F_1\in \Theta^{\prime}_\varepsilon(F_{0})}\sqrt{k_{n}} |R_{k_{n},n} - 1| \overset{d}{\longrightarrow} +\infty.
\]
\end{enumerate}
\end{theorem}

\section{Auxiliary results and proofs.}

\subsection{Auxiliary results.}

Since $R_{n}$ depends on the higher order statistics we cannot immediately use
independence of the random variables $(X_{1},\ldots,X_{n}).$ Therefore
consider the conditional distribution of $R_{n}$ given $X_{(n-k)}=q$ applying
the following lemma.

\begin{lemma}
\label{L1} (\cite{Dehaan}) Let $X,X_{1},\ldots,X_{n}$ be i.i.d. random
variables with common distribution function $F,$ and let $X_{(1)}\leq
\ldots\leq X_{(n)}$ be the $n$th order statistics. For any $k=1\ldots n-1$,
the conditional joint distribution of $\{X_{(i)}\}_{i=n-k+1}^{n}$ given
$X_{(n-k)}=q$ is equal to the (unconditional) joint distribution of the
corresponding set $\{X_{(i)}^{\ast}\}_{i=1}^{k}$ of order statistics for
i.i.d. random variables $\{X_{i}^{\ast}\}_{i=1}^{k}$ having the distribution function
\end{lemma}

\[
F_{q}(x)=P(X\leq x|X>q)=\frac{F(x)-F(q)}{1-F(q)},\ \ \ \ x>q.
\]
We call $F_{q}(x),$ $x>q,$ the tail distribution function linked with the
distribution function $F.$ Consider two continuous distribution functions $F$
and $G$ and a random variable $\xi_{q}$ with distribution function $G_{q},$
where $q\in\mathbb{R}$ is some parameter. Let
\[
\eta_{q}=\ln\left(  \frac{1-F(q)}{1-F(\xi_{q})}\right)  .
\]
Clear, $\eta_{q}\geq0$ for all $q\in\mathbb{R}.$ \newline The crucial point in
the proof of Theorem \ref{Th2} is studying of asymptotical behavior of
$\eta_{q}.$

\begin{proposition}
\label{P1} Let $F_{q}$ and $G_{q}$ are tail distribution functions of $F$ and
$G$ respectively. Then

\begin{description}
\item[(i)] If for some $x_{0}$, $q>x_{0},$ and any $x>q$, $F_{q}(x)=G_{q}(x),$
then  $\eta_{q}$ is standard exponential.

\item[(ii)] $G_{q}(x)\geq F_{q}(x)$ for any $x>q$ if and only if $\eta_{q}$ is
stochastically smaller than a standard exponential random variable.\newline%
$G_{q}(x)\leq F_{q}(x)$ for any $x>q$ if and only if $\eta_{q}$ is
stochastically larger than a standard exponential random variable.

\item[(iii)] $G_{q}(x)\geq F_{q}(x)$ for any $x>q\geq x_{0}$ and some $x_{0}$
if and only if $(1-G(x))/(1-F(x))$ is nonincreasing function as $x>x_{0}.$
\end{description}
\end{proposition}

\subsection{Proof of Proposition \ref{P1}.}

\emph{(i)} Let $F_{q}(x)=G_{q}(x)$ for all $x>q,$ then we have for the
distribution function of $\eta_{q}$,
\[
P(\eta_{q}\leq y)=P\left(  \ln\left(  \frac{1-F(q)}{1-F(\xi_{q})}\right)  \leq
y\right)  =P\left(  \frac{1-F(q)}{1-F(\xi_{q})}\leq e^{y}\right)  =
\]%
\begin{equation}
=P\left(  F(\xi_{q})\leq1-(1-F(q))e^{-y}\right)  =P\left(  \xi_{q}\leq
F^{\leftarrow}\left(  1-\frac{1-F(q)}{e^{y}}\right)  \right)  .\label{Petaq}%
\end{equation}
Furthermore, for the same $x$,
\[
P\left(  \xi_{q}\leq F^{\leftarrow}\left(  1-\frac{1-F(q)}{e^{y}}\right)
\right)  =\frac{F\left(  F^{\leftarrow}\left(  1-\frac{1-F(q)}{e^{y}}\right)
\right)  -F(q)}{1-F(q)}=1-e^{-y}.
\]
\emph{(ii)} Now assume that for all $x>q$ and some $q\in\mathbb{R}$,
$G_{q}(x)\geq F_{q}(x).$ Then from \eqref{Petaq}, since $1-(1-F(q))e^{-y}\geq
F(q)$ for all  $y\geq0$ it follows that
\begin{equation*}
P(\eta_{q}\leq y)=\frac{G\left(  F^{\leftarrow}\left(  1-\frac{1-F(q)}{e^{y}%
}\right)  \right)  -G(q)}{1-G(q)}\geq\end{equation*}
\begin{equation}\frac{F\left(  F^{\leftarrow}\left(
1-\frac{1-F(q)}{e^{y}}\right)  \right)  -F(q)}{1-F(q)}=1-e^{-y}.\label{F<G}%
\end{equation}
Conversely, assume that $\eta_{q}$ is stochastically smaller than a standard
exponential random variable, that is, $P(\eta_{q}\leq y)\geq1-e^{-y}$ for all
$y\geq0.$ With \eqref{Petaq} we get that
\[
\frac{G\left(  F^{\leftarrow}\left(  1-\frac{1-F(q)}{e^{y}}\right)  \right)
-G(q)}{1-G(q)}\geq1-e^{-y}\Longleftrightarrow\frac{1-G\left(  F^{\leftarrow
}\left(  1-\frac{1-F(q)}{e^{y}}\right)  \right)  }{1-G(q)}\leq e^{-y}%
\]%
\[
\Longleftrightarrow G\left(  F^{\leftarrow}\left(  1-\frac{1-F(q)}{e^{y}}\right)  \right)
\leq1-\frac{1-G(q)}{e^{y}}\Longleftrightarrow\]
\[ F^{\leftarrow}\left(
1-\frac{1-F(q)}{e^{y}}\right)  \leq G^{\leftarrow}\left(  1-\frac
{1-G(q)}{e^{y}}\right)  .
\]
Denote $z_{F}=F^{\leftarrow}\left(  1-e^{-y}(1-F(q))\right)  $ and
$z_{G}=G^{\leftarrow}\left(  1-e^{-y}(1-G(q))\right)  .$ Since $F(z_{F}%
)=1-e^{-y}(1-F(q))$ and $G(z_{G})=1-e^{-y}(1-G(q)),$ we have,
\[
e^{-y}=\frac{1-G(z_{G})}{1-G(q)}=\frac{1-F(z_{F})}{1-F(q)}.
\]
Further, since $z_{F}\leq z_{G}$ then
\[
\frac{1-F(z_{F})}{1-F(q)}=\frac{1-G(z_{G})}{1-G(q)}\leq\frac{1-G(z_{F}%
)}{1-G(q)}.
\]
This observation completes the proof since $z_{F}\in\lbrack q,\infty).$ The
proof of the second assertion is similar.\newline\newline\emph{(iii)} We
have,
\[
\frac{G(x)-G(q)}{1-G(q)}\geq\frac{F(x)-F(q)}{1-F(q)}\ \ \forall x>q\geq
x_{0}\Longleftrightarrow\frac{1-G(x)}{1-G(q)}\leq\frac{1-F(x)}{1-F(q)}%
\ \ \forall x>q\geq x_{0}\Longleftrightarrow
\]%
\[
\frac{1-G(x)}{1-F(x)}\leq\frac{1-G(q)}{1-F(q)}\ \ \forall x>q\geq
x_{0}\Longleftrightarrow\frac{1-G(x)}{1-F(x)}\text{ is nonincreasing for all
}x>x_{0}.\ \ \blacksquare
\]
%It means that $\eta_q$ is stochastically smaller $\forall q>x_0$ than a standard exponential random variable.\\
%So $\eta_q$ is stochastically larger $\forall q>x_0$ than a standard exponential random variable
%and if $E\eta_q<+\infty$ for some $q>x_0$ then $E\eta_q>1.$
%\emph{(iv)} Suppose \eqref{conv} hold then
%\[\int_q^{\infty} \ln (1- F(x)) dG(x)<+\infty\]
%as $q\to\infty.$ See, this integral exists for some $q$ if and only if $E \eta_q<+\infty.$
%Proof of \eqref{conv2} is the same.

\subsection{Proof of Theorem \ref{Th1}.}

Under the conditions of Theorem \ref{Th1}, $F_{0}(X_{1})$ is uniformly
distributed on $[0,1]$, that is, $F_{0}(X_{1})$ $\sim U[0,1],$ hence
$-\ln(1-F_{0}(X))$ is standard exponential random variable. It follows from
R\'{e}nyi's representation (see \cite{Dehaan}), that
\[
\left\{  -\ln(1-F_{0}(X_{(n-i)}))+\ln(1-F_{0}(X_{(n-k)}))\right\}
_{i=0}^{k-1}\overset{d}{=}\left\{  \sum\limits_{j=i+1}^{k}\frac{E_{n-j+1}}%
{j}\right\}  _{i=0}^{k-1},
\]
where $E_{1},E_{2}\ldots$ are independent standard exponential variables.
Therefore the distribution of the left-hand side does not depend on $n$ and
\[
\left\{  -\ln(1-F_{0}(X_{(n-i)}))+\ln(1-F_{0}(X_{(n-k)}))\right\}
_{i=0}^{k-1}\overset{d}{=}\left\{  E_{(k-i)}\right\}  _{i=0}^{k-1},
\]
where $E_{(1)}\leq\ldots\leq E_{(k)}$ are the $n$th order statistics of the
sample $\{E_{i}\}_{i=1}^{k}.$ Finally we have,
\[
\sqrt{k}(R_{k,n}-1)\overset{d}{=}\sqrt{k}\left(  \frac{1}{k}\sum
\limits_{i=0}^{k-1}E_{(k-i)}-1\right)  =\sqrt{k}\left(  \frac{1}{k}%
\sum\limits_{j=1}^{k}E_{j}-1\right)  ,
\]
and the assertion follows from the Central Limit Theorem.

\subsection{Proof of Theorem \ref{Th2}.}
We first prove (i).
The steps of the proof are similar to corresponding steps in \cite{Rodionov1}
and \cite{Rodionov2}. Consider asymptotic behavior of $R_{k_{n},n}$ as
$n\rightarrow\infty.$ Denote
\[
Y_{i}=\ln(1-F_{0}(q))-\ln(1-F_{0}(X_{i}^{\ast})),
\]
where $\{X_{i}^{\ast}\}_{i=1}^{k_{n}}$ are i.i.d. random variables introduced
in Lemma \ref{L1} with the distribution function
\[
F_{q}(x)=\frac{F_{1}(x)-F_{1}(q)}{1-F_{1}(q)},\ \ q<x.
\]
Taking $F=F_{0}$ and $G=F_{1}$ we have, $Y_{i}\overset{d}{=}\eta_{q},$
$i\in\{1,\ldots,k_{n}\}$.  Notice that, in view of Lemma 1, the joint
distribution of order statistics $\{Y_{(i)}\}_{i=1}^{k_{n}}$ of the sample
$\{Y_{j}\}_{i=1}^{k_{n}}$ is equal to the joint conditional distribution of
order statistics $\{Z_{(j)}\}_{j=1}^{k_{n}}$ of  $\{Z_{j}\}_{j=1}^{k_{n}}$
given $X_{(n-k_{n})}=q,$ where
\[
Z_{j}=\ln(1-F_{0}(X_{(n-k_{n})}))-\ln(1-F_{0}(X_{(n-j+1)})),\ j=1,...,k_{n}.
\]
Clear,
\[
R_{n,k_{n}}=\frac{1}{k_{n}}\sum_{i=1}^{k_{n}}Z_{i}.
\]
So, the conditional distribution of $R_{k_{n},n}$ given $X_{(n-k)}=q$ is equal
to the distribution of $\frac{1}{k_{n}}\sum_{i=1}^{k_{n}}Y_{i}.$ Further,
distribution functions $F_{1}$ and $F_{0}$ satisfy $B(F_{0},F_{1})$ or
$B(F_{1},F_{0}).$ First suppose that the condition $B(F_{0},F_{1})$ holds for
some $\varepsilon>0$ and $x_{0}.$ Since $x^{\ast}=+\infty,$ $X_{(n-k_{n}%
)}\rightarrow+\infty$ a.s., we may consider the case $q>x_{0}$ only.
Proposition \ref{P1} (iii) implies, that
\[
\frac{1-F_{1}(x)}{1-F_{1}(x_{0})}\geq\frac{(1-F_{0}(x))^{1-\varepsilon}%
}{(1-F_{0}(x_{0}))^{1-\varepsilon}},\ \ x>x_{0}.
\]
With \eqref{F<G}, we get that,
\[
P(Y_{1}\leq x)=1-\frac{1-F_{1}\left(  F_{0}^{\leftarrow}\left(  1-\frac
{1-F_{0}(q)}{e^{x}}\right)  \right)  }{1-F_{1}(q)}\leq\]
\[1-\frac{\left(
1-F_{0}\left(  F_{0}^{\leftarrow}\left(  1-\frac{1-F_{0}(q)}{e^{x}}\right)
\right)  \right)  ^{1-\varepsilon}}{(1-F_{0}(q))^{1-\varepsilon}%
}=1-e^{-(1-\varepsilon)x},
\]
hence $Y_{1}$ is stochastically larger than a random variable $E\sim
Exp(1-\varepsilon),$ write $Y_{1}\gg E.$ Further, let $E_{1},\ldots,E_{k_{n}}$
are i.i.d. random variables with distribution function
$H(x)=1-e^{-(1-\varepsilon)x},$ then
\begin{equation}
\sqrt{k_{n}}\left(  \frac{1}{k_{n}}\sum_{i=1}^{k_{n}}Y_{i}-1\right)  \gg
\sqrt{k_{n}}\left(  \frac{1}{k_{n}}\sum_{i=1}^{k_{n}}E_{i}-1\right)
.\label{ll}%
\end{equation}
Since \eqref{ll} holds for all $q>x_{0},$ and $X_{(n-k_{n})}\rightarrow
+\infty$ a.s. as $n\to\infty$, we have under the conditions of Theorem \ref{Th2}, that
\begin{equation}
\sqrt{k_{n}}(R_{k_{n},n}-1)\gg\sqrt{k_{n}}\left(  \frac{1}{k_{n}}\sum
_{i=1}^{k_{n}}E_{i}-1\right)  .\label{LL}%
\end{equation}
It follows from Lindeberg-Feller theorem, that
\[
(1-\varepsilon)\sqrt{k_{n}}\left(  \frac{1}{k_{n}}\sum_{i=1}^{k_{n}}%
E_{i}-\frac{1}{1-\varepsilon}\right)  \xrightarrow{d}\xi\sim
N(0,1),\ \ n\rightarrow\infty,
\]
therefore
\begin{equation}
\sqrt{k_{n}}\left(  \frac{1}{k_{n}}\sum_{i=1}^{k_{n}}E_{i}-1\right)
\xrightarrow{P}+\infty,\ \ n\rightarrow\infty.\label{ref2}
\end{equation}
Finally, with \eqref{LL}, we have,
\[
\sqrt{k_{n}}(R_{k_{n},n}-1)\xrightarrow{P}+\infty,\ \ n\rightarrow\infty.
\]
If the condition $B(F_{0},F_{1})$ holds, then
\[
\sqrt{k_{n}}(R_{k_{n},n}-1)\xrightarrow{P}-\infty,\ \ n\rightarrow\infty,
\]
and the proof is similar. The second assertion easily follows from \eqref{LL} and \eqref{ref2}.

\subsection{Proof of Theorem \ref{Th3}.}
Firstly we prove (i).
Denote $\overline{F}(x)=1-F(x).$ In notation of the proof of Theorem
\ref{Th2}, find the distribution of $Y_{1}.$ First assume that $C(F_{0}%
,F_{1})$ holds. With \eqref{F<G} and Proposition \ref{P1} (iii) we have,
\[
P(Y_{1}\leq x)=1-\frac{\overline{F_{1}}\left(  \overline{F_{0}}^{\leftarrow
}(\overline{F_{0}}(q){e^{-x}})\right)  }{\overline{F_{1}}(q)}\leq\]
\[
1-\frac{\overline{F_{0}}(q)e^{-x}\left(  -\ln(\overline{F_{0}}(q)e^{-x}%
)\right)  ^{\varepsilon}}{\overline{F}(q)(-\ln\overline{F_{0}}%
(q))^{\varepsilon}}=1-e^{-x}\left(  1+\frac{x}{-\ln\overline{F}_{0}%
(q)}\right)  ^{\varepsilon}.
\]
For $\varepsilon,c\in(0,1),$
\[
(1+cx)^{\varepsilon}\geq1+c\varepsilon-c\varepsilon e^{-x},\ \ x\geq0,
\]
and $G(x)=1-e^{-x}(1+c\varepsilon-c\varepsilon e^{-x})$ is the distribution
function. Hence,
\[
P(Y_{1}\leq x)\geq1-e^{-x}-\frac{\varepsilon}{-\ln\overline{F_{0}}%
(q)}(1-e^{-x}).
\]
Further, let $\zeta,\zeta_{1},\ldots,\zeta_{k_{n}}$ be i.i.d. random variables
with this distribution function. Therefore, like the proof of Theorem
\ref{Th2},
\begin{equation}
\sqrt{k_{n}}\left(  \frac{1}{k_{n}}\sum_{i=1}^{k_{n}}Y_{i}-1\right)  \gg
\sqrt{k_{n}}\left(  \frac{1}{k_{n}}\sum_{i=1}^{k_{n}}\zeta_{i}-1\right)
.\label{ll1}%
\end{equation}
Clear,
\[
E\zeta=1+\frac{\varepsilon}{-2\ln\overline{F_{0}}(q)},\ \ Var\,\zeta=1-\left(
\frac{\varepsilon}{-2\ln\overline{F_{0}}(q)}\right)  ^{2},
\]
so we have,
\begin{equation} \label{ref3}
\sqrt{k_{n}}\left(  \frac{1}{k_{n}}\sum_{i=1}^{k_{n}}\zeta_{i}-1\right)
=\sqrt{k_{n}}\left(  \frac{1}{k_{n}}\sum_{i=1}^{k_{n}}\zeta_{i}-E\zeta\right)
+\sqrt{k_{n}}\frac{\varepsilon}{-2\ln\overline{F_{0}}(q)}.
\end{equation}
Consider now the statistic $\sqrt{k_{n}}/\ln\overline{F_{0}}(X_{(n-k_{n})}).$
Denote $R_{i}=\overline{F_{1}}(X_{i}),$ $i=1,\ldots,n.$ Since $F_{1}$ is
continuous, $R_{1},\ldots,R_{n}$ are i.i.d. standard uniform random variables
and $R_{(k_{n})}=\overline{F_{1}}(X_{(n-k_{n})}).$ Theorem 2.2.1 \cite{Dehaan}
implies, that
\begin{equation}
\frac{n}{\sqrt{k_{n}}}\left(  R_{(k_{n})}-\frac{k_{n}}{n}\right)
\xrightarrow{d}N(0,1),\ n\to\infty. \label{last}%
\end{equation}
Using the delta method (see \cite{Resnick}) for the function $f(x)=-x/\ln x,$
we have
\[
\frac{n}{\sqrt{k_{n}}}\left(  \frac{R_{(k_{n})}}{-\ln R_{(k_{n})}}-\frac
{k_{n}/n}{-\ln(k_{n}/n)}\right)  \xrightarrow{P}0,\ n\to\infty,
\]
since under the conditions of theorem
\[
f^{\prime}\left(  \frac{k_{n}}{n}\right)  =-\frac{1}{\ln(n/k_{n})}+\frac
{1}{(\ln(n/k_{n}))^{2}}\rightarrow0,\ n\rightarrow\infty.
\]
Further,
\[
\frac{n}{\sqrt{k_{n}}}\left(  \frac{R_{(k_{n})}}{\ln R_{(k_{n})}}-\frac
{k_{n}/n}{\ln(k_{n}/n)}\right)  =\]
\[\frac{n}{\sqrt{k_{n}}}\left(  \frac
{R_{(k_{n})}}{\ln R_{(k_{n})}}-\frac{k_{n}/n}{\ln(R_{(k_{n})})}\right)
+\sqrt{k_{n}}\left(  \frac{1}{\ln R_{(k_{n})}}-\frac{1}{\ln(k_{n}/n)}\right)
,
\]
and \eqref{last} implies that the first summand in the right hand side tends
to $0$ in probability. Therefore,
\[
\sqrt{k_{n}}\left(  \frac{1}{\ln R_{(k_{n})}}-\frac{1}{\ln(k_{n}/n)}\right)
\xrightarrow{P}0,\ n\to\infty,
\]
and under the conditions of Theorem \ref{Th3},
\[
\frac{\sqrt{k_{n}}}{-\ln R_{(k_{n})}}=\sqrt{k_{n}}\left(  \frac{1}{-\ln
R_{(k_{n})}}-\frac{1}{-\ln(k_{n}/n)}\right)  +\frac{\sqrt{k_{n}}}{-\ln
(k_{n}/n)}\xrightarrow{P}+\infty,\ n\to\infty.
\]
On the other hand, from \eqref{new} it  follows that
\begin{equation} \label{ref4}
\frac{\sqrt{k_{n}}}{-\ln\overline{F_{0}}(X_{(n-k_{n})})}=\frac{\sqrt{k_{n}}%
}{-\ln\overline{F_{0}}\left(  \overline{F_{1}}^{\leftarrow}(R_{(k_{n}%
)})\right)  }\geq\frac{\sqrt{k_{n}}}{-\delta^{-1}\ln R_{(k_{n})}%
}\xrightarrow{P}+\infty
\end{equation}
as $n\rightarrow\infty.$ Further, it follows from the Law of large numbers for
triangular arrays (see \cite{Mik}), that for any $\epsilon>0$
\[
\sqrt{k_{n}}\left(  \frac{1}{k_{n}}\sum_{i=1}^{k_{n}}\zeta_{i}-1\right)
=o_{P}(k_{n}^{\epsilon}),\ \ n\rightarrow\infty.
\]
It means that the term in the left hand side is asymptotically smaller in
probability than $k_{n}^{\epsilon}.$ Hence for any $q,$ given $X_{(n-k_{n}%
)}=q$
\[
\sqrt{k_{n}}\left(  \frac{1}{k_{n}}\sum_{i=1}^{k_{n}}Y_{i}-1\right)
\xrightarrow{P}+\infty,\ n\to\infty,
\]
and finally,
\[
\sqrt{k_{n}}\left(  R_{k_{n},n}-1\right)  \xrightarrow{P}+\infty,\ n\to\infty.
\]
If the condition $C(F_{1},F_{0})$ holds, then
\[
\sqrt{k_{n}}(R_{k_{n},n}-1)\xrightarrow{P}-\infty,\ n\rightarrow\infty,
\]
and the proof is the same. The second assertion clearly follows from \eqref{ll1}, \eqref{ref3} and \eqref{ref4}.

\end{document}